\documentclass{article}
\usepackage{spconf,amsmath,amssymb,amsthm,graphicx}
\usepackage[ruled, linesnumbered]{algorithm2e}
\usepackage[dvipsnames]{xcolor}
\usepackage{hyperref}
\usepackage{cleveref}
\usepackage{subcaption}

\usepackage{bbm}

\usepackage{tikz}
\usetikzlibrary{shapes}

\newtheorem{theorem}{Theorem}
\newtheorem{lemma}{Lemma}
\newtheorem{proposition}{Proposition}

\newif\iflong
\longtrue
\newcommand{\inExtendedVersion}[1]{\iflong #1\fi}
\newcommand{\onlyInShortVersion}[1]{\iflong \else #1 \fi}

\title{Blind extraction of equitable partitions from graph signals}
\name{Michael Scholkemper, Michael T. Schaub\thanks{We acknowledge partial funding from the Ministry of Culture and Science of North Rhine-Westphalia (NRW Rückkehrprogramm) and the Excellence Strategy of the Federal Government and the Länder. \inExtendedVersion{Copyright 2022 IEEE. Personal use of this material is permitted. Permission from IEEE must be obtained for all other uses, in any current or future media, including reprinting/republishing this material for advertising or promotional purposes, creating new collective works, for resale or redistribution to servers or lists, or reuse of any copyrighted component of this work in other works.}}}
\address{Department of Computer Science, RWTH Aachen University, Germany}
\begin{document}
\ninept
\maketitle

\begin{abstract}

Finding equitable partitions is closely related to the extraction of graph symmetries and of interest in a variety of applications context such as node role detection, cluster synchronization, consensus dynamics, and network control problems.
In this work we study a blind identification problem in which we aim to recover an equitable partition of a network without the knowledge of the network's edges but based solely on the observations of the outputs of an unknown graph filter.
Specifically, we consider two settings.
First, we consider a scenario in which we can control the input to the graph filter and present a method to extract the partition inspired by the well known Weisfeiler-Lehman (color refinement) algorithm.
Second, we generalize this idea to a setting where only observe the outputs to random, low-rank excitations of the graph filter, and present a simple spectral algorithm to extract the relevant equitable partitions.
Finally, we establish theoretical bounds on the error that this spectral detection scheme incurs and perform numerical experiments that illustrate our theoretical results and compare both algorithms.
\end{abstract}
\begin{keywords}
Equitable partitions, Weisfeiler Lehman algorithm, spectral analysis, topology inference, graph symmetry
\end{keywords}
\section{Introduction}
\label{sec:intro}

Networks have become a powerful abstraction for complex systems~\cite{strogatz2001exploring,newman2018networks}. 
To comprehend such networks, we often seek patterns in their connections, which would allow us to comprehend the system in simpler terms.
A common theme is to divide the nodes---and by extension the units of the underlying system---into groups of similar nodes.
For instance, in the context of community detection, we consider nodes as similar if they are tightly-knit together, or share similar neighborhoods~\cite{fortunato2010community}.
This notion of node similarity is thus bound to the specific position of the nodes in the graph, i.e., the identity of their neighboring nodes.
In contrast, we may want to split nodes into groups according to whether they play a similar \emph{role} in the graph~\cite{rossi2020proximity}, irrespective of their exact position.
As an example, consider a division into hubs and peripheral nodes according to their degree, a split for which the exact identity of the neighboring nodes is not essential.
While in this specific example defining a degree-similarity measure between nodes is simple, how to define a similarity measure between nodes in a position independent manner is a non-trivial question in general.

Rather than trying to identify similar nodes, many traditional approaches in social network analysis consider the definition of nodes roles based on \emph{exact node equivalences}, such as regular equivalence or automorphic equivalence~\cite{brandes2005network}.
A specific form of such a partition into exact node equivalence classes is an \emph{equitable partition} (EP), which may be intuitively defined recursively as sets of nodes that are connected to the same number of equivalent nodes.
These EPs generalize orbit partitions related to automorphic equivalence and are thus closely related to graph symmetries.
Knowledge of EPs can thus, e.g., facilitate the computation of network statistics such as centrality measures~\cite{sanchez2020exploiting}.
As they are associated to certain spectral signatures, EPs are also relevant for the study of dynamical processes on networks such as cluster synchronization~\cite{pecora2014cluster,schaub2016graph}, consensus dynamics~\cite{yuan2013decentralised}, and network control problems~\cite{martini2010controllability}.

Motivated by this interplay between network dynamics and EPs, in this work we ask the question: \emph{Can we detect the presence of an EP in a network solely based on a small number of nodal observations of a dynamical process acting on the network?}
For this, we adopt a graph signal processing perspective \cite{ortega2018graph}, in which we model the dynamics as a graph signal filtered by an unknown filter representing the dynamics.
Our task is then to recover the EP solely based on a small number of outputs of this filter.


\noindent\textbf{Related work.}
Network topology inference has been studied extensively in the literature~\cite{dong2016learning,kalofolias2016learn}.
Inferring the complete topology of a network can however require a large number of samples and may thus be infeasible in practice.
A relatively recent development is the idea to bypass this inference of the exact network topology and directly estimate network characteristics in the form of community structure~\cite{wai2018community,roddenberry2020exact,schaub2020blind} or centrality measures \cite{roddenberry2020blind,he2020estimating,he2021detecting} from graph signals.
Learning graph characteristics directly in this way benefits from a better sample complexity, since only a low-dimensional set of spectral features must be inferred rather than the whole graph.
This manuscript falls squarely within this line of work, but focuses on EPs as a different network feature, which results in some distinct challenges.
Specifically, we cannot rely on the estimation of a dominant invariant subspace, but must estimate and select a subset of relevant eigenvectors from the whole spectrum of the graph.

\noindent\textbf{Contributions and outline.} 
We present two algorithms to tackle the problem of extracting the coarsest EP from the observation of graph signals under two different scenarios.
First, we consider a scenario where we can control the input to the (unknown) graph filter, while having no access to the graph.
We present an algorithm that exactly recovers the coarsest EP in this setting.
Second, we consider a fully ``blind'' estimation problem, where we only have access to noisy, random, low-rank excitations of the graph filter. 
For this we derive a simple spectral algorithm and derive theoretical error bounds for its performance under certain assumptions.
Finally, we illustrate our results and compare the two algorithms.

\section{Notation}\label{sec:prelim}


\textbf{Graphs.} A simple graph $G = (V,E)$ consists of a set of nodes $V$ and a set of edges $E = \{uv \mid u,v \in V\}$.
The neighborhood $N(v) = \{x \mid vx \in E\}$ of a node $v$ is the set of all nodes connected to $v$.
A graph is \emph{undirected} if $uv \in E \iff vu \in E$.
We may consider more general (non-simple) graphs that allow for directed edges, self loops $vv \in E$, or assign positive weights $w:E \rightarrow \mathbb{R}_+$ to the edges of the graph, rendering it a \emph{weighted graph}.

\noindent\textbf{Matrices.}
For a matrix $M$, $M_{i,j}$ is the component in the $i$-th row and $j$-th column. 
We use $M_{i,\_}$ to denote the $i$-th row vector of $M$ and $M_{\_,j}$ to denote the $j$-th column vector.
The span of a matrix $M\in\mathbb{R}^{n\times m}$ is defined as the set $\operatorname{span}(M) = \{x \in \mathbb{R}^{n} \mid x = M v \text{ for some } v\in \mathbb{R}^m\}$. 
For a matrix $X$, we denote by $X \in \operatorname{span}(M)$ if $X_{\_,j}\in\operatorname{span}(M)$ for all $j$.
$\mathbb{I}_n$ is the identity matrix and $\mathbbm{1}_n$ the all-ones vector of size $n$ respectively.
Given a graph $G = (V,E)$, we identify the set of nodes with $\{1, \ldots, n\}$. 
An \emph{adjacency matrix} of a given graph is a matrix $A$ with entries $A_{u,v} = 0$ if $uv \notin E$ and $A_{u,v} = w(uv)$ otherwise, where we set $w(uv)=1$ for unweighted graphs for all $uv\in E$.

\noindent\textbf{Equitable partitions.}
A \emph{partition} $\mathcal{C} = \{C_1, ..., C_k\}$ of a graph into $k$ classes splits the nodes into disjoint sets $C_i$ such that $C_i \cap C_j = \emptyset$ for $i \neq j$ and $V = \bigcup_i C_i$.
An \emph{equitable partition} (EP) is a partition such that within each class $C_i$ the neighborhood of each node is partitioned into parts of the same size for all classes $C_j$.
Formally, in an EP it holds for all nodes $v, u \in C_i$ in the same class that:
\begin{equation}
    |\{N(v) \cap C_j\}| = |\{N(u) \cap C_j\}| \text{ for all } j \in \{1, ..., k\}.
\end{equation}
A standard algorithm to detect EPs is the Weisfeiler-Lehman algorithm~\cite{weisfeiler1968reduction}, a combinatorial algorithm which iteratively refines a coloring of the nodes until a stable coloring is reached. It finds the \textit{coarsest EP} (cEP), meaning that with the fewest classes. 
Note that it is typically this cEP one aims to find, as it provides the largest reduction in complexity for describing the structure of the graph. Furthermore, any graph will have a trivial finest EP with $|V|$ classes. 


\section{Equitable partitions and eigenvectors}
\label{sec:counterexample}
\begin{figure}[t]
	\small
     	\begin{subfigure}[b]{0.23\textwidth}
         \centering
         \begin{tikzpicture}[->,>=stealth,shorten >=1pt,auto,node distance=1cm,
			thick,main node/.style={circle,draw,minimum size=0.3cm,inner sep=0pt}]
			
			\node[main node, cyan] (5) {$5$};
			\node[main node, cyan] (6) [above right of=5]  {$6$};
			\node[main node, orange] (1) [below right of=5] {$1$};
			\node[main node, orange] (2) [right of=6] {$2$};
			\node[main node, orange] (3) [below right of=2] {$3$};
			\node[main node, orange] (4) [below left of=3] {$4$};

			\path[-]
			(1) edge node {2} (5)
			(5) edge node {4} (6)
			edge node {} (1)
			(6) edge node {2} (2)
		  edge node {} (5)
			(2) edge node {4} (3)
			edge node {} (6)
			(3) edge node {4} (4)
			edge node {} (2)
			(4) edge node {} (3)
			(1) edge [loop below] node {4} (1)
	
			(4) edge [loop below] node {4} (4)
			(5) edge [loop below] node {3} (5)
			(6) edge [loop below] node {3} (6)
	
			;
			\end{tikzpicture}
         \caption{}
         \label{fig:y equals x}
     	\end{subfigure}
     	\hfill
     	\begin{subfigure}[b]{0.23\textwidth}
			\centering
			\begin{tikzpicture}[->,>=stealth,shorten >=1pt,auto,node distance=1cm,
			   thick,main node/.style={circle,draw,minimum size=0.3cm,inner sep=0pt}]
			   
			   \node[main node, green] (1) {$1$};
			   \node[main node, green] (2) [below left of=1]  {$2$};
			   \node[main node, orange] (3) [below right of=1] {$3$};
			   \node[main node, gray] (4) [below right of=2] {$4$};

			   \path[-]
			   (1) edge node {} (2)
			   edge node {} (3)
			   (2) edge node {} (1)
			   edge node {} (3)
			   edge node {} (4)
			   (3) edge node {} (1)
			   edge node {} (2)
			   edge node {} (4)
			   (4) edge node {} (2)
			   edge node {} (3)
			   (3) edge [loop below] node {2} (1)
	   
			   (4) edge [loop below] node {2} (4)
	   
			   ;
			   \end{tikzpicture}
			\caption{}
			\label{fig:y equals x2}
			\end{subfigure}
			\hfill

	\caption{Examples of graphs where the structural eigenvectors of the cEP behave irregularly. The depicted graphs are loopy to make the examples small. However, we can think of this as the condensed graphs $A^\pi$ of some larger simple graph. (a) $[1,1,1,1,2,2]$ is the perron vector of the graph. However, the cEP has $6$ classes. (b) This graph also only has one EP, but the adjacency matrix is singular. Therefore, an eigenvector with eigenvalue $0$ is part of the cEP.}
	\label{fig:counterexample}
\end{figure}
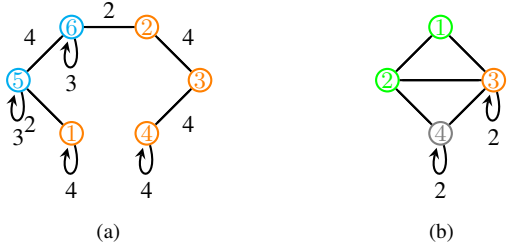

In this section we recap how the presence of an equitable partition manifests itself in the spectral properties of the adjacency matrix.
We will subsequently use these spectral properties to detect the presence of a cEP from the output of a graph filter.


Throughout the paper, we consider an undirected graph $G$ with a cEP that has $k$ classes encoded by the indicator matrix $H$ with entries $H_{ij} = 1$ if $i\in C_j$ and $H_{ij}= 0$ otherwise. 
It now holds that
\begin{equation}\label{eq:defEP}
    AH = H(H^\top H)^{-1} H ^\top A H =: HA^\pi 
\end{equation}
where $A^\pi$ is the adjacency matrix of the so-called \textit{quotient graph} associated to the EP.
Note that $A^\pi$ is not necessarily simple nor undirected.
The converse of \cref{eq:defEP} also holds. 
If there exists an indicator matrix $H \in \{0,1\}^{n \times k}$ as defined above and $AH = HA^\pi$, then the graph has an EP with $k$ classes as indicated by $H$.

\noindent\textbf{Spectral signatures of EPs}
The algebraic characterization of EPs in terms of~\Cref{eq:defEP} has some noteworthy consequences.
First, the eigenvalues of $A^\pi$ are a subset of the eigenvalues of $A$.
Second, the eigenvectors of $A^\pi$ can be scaled up by $H$ to become eigenvectors of $A$.
Both of these statements are implied by the following argument: Let $A^\pi v = \lambda v$ define an eigenpair of $A^\pi$, then $AHv = HA^\pi v = \lambda Hv$ is an eigenpair of $A$.

Thus, $A$ has $k$ (unnormalized) eigenvectors of the form $z=Hv$ that are by construction block constant on the classes of the EP.\@
We call these eigenvectors of $A$ associated with the cEP \textit{structural eigenvectors}.
As $A$ is symmetric, these eigenvectors form an orthogonal basis of a closed subspace. 
Moreover, as the structural eigenvectors are block-wise constant, the same subspace is also spanned by the indicator vectors of the classes, i.e., the columns of~$H$.

This motivates the following intuitive idea to find the cEP \emph{spectrally} by analyzing the eigenvectors of the adjacency matrix $A$:
\textit{Find the smallest number $k$ of eigenvectors that are block constant on the same $k$ blocks}.
The vectors indicating the blocks of these eigenvectors will then correspond to the columns of the indicator matrix $H$ for the cEP of the graph.

Since any structural eigenvector is of the form $z=Hv$, we may even hope to find the cEP by simply computing a single structural eigenvector, provided that the entries of the eigenvector $v$ of the quotient graph are all distinct. 
Given such an eigenvector $z$, we could then simply read off the node-equivalences classes by checking whether or not two nodes are assigned the same value in the vector $z$.
In fact, we can always identify at least one structural eigenvector of the cEP easily using the following proposition, which follows from the non-negativity of the adjacency and the algebraic characterization of the cEP given above.
\begin{proposition}
	\label{prop:1}
    If the Perron vector (dominant eigenvector) of a graph is unique, then
	it is a structural eigenvector.
\end{proposition}

However, there are several obstacles to implement the above algorithmic idea.
First, even if we are given a structural eigenvector of the cEP, simply grouping the nodes according to their entries in the eigenvector may not reveal the cEP.
For instance, focussing only on the Perron vector may not yield the correct cEP, as nodes in different classes may still be assigned the same value. 
Second, the eigenvectors associated to the cEP may be associated to any eigenvalue of $A$, i.e., they are not necessarily the dominant eigenvectors of $A$.

Both of these problematic cases are shown in~\Cref{fig:counterexample}.
Figure~\Cref{fig:counterexample}a provides an example where the dominant eigenvector (Perron vector) has fewer than $k$ distinct entries.
Figure~\Cref{fig:counterexample}b shows an eigenvector of the cEP with small eigenvalues (in this case zero).
While there area thus always exactly $k$ structural eigenvectors of the cEP, we must first determine which eigenvectors are indeed structural eigenvectors of the cEP and which are not.

\section{Scenario I: Blind but in control}\label{sec:simWL}
In the following, we consider a setup in which we cannot see the graph directly, but can sample the input/output behavior of a graph filter in form of a matrix polynomial $f(A)= \sum_k h_k A^k$, where $h_k\in\mathbb{R}$ are filter coefficients.
To illustrate our main ideas, we will restrict ourselves for this section to the simplest case, where $f(A)=A$ and we can control the input to the filter.
In~\Cref{sec:main} we will then concern ourselves with a fully ``blind'' cEP inference problem, where we have no control over the inputs.

For now let us assume we observe the outputs $y=Ax\in\mathbb{R}^n$ to a set of inputs $x$ we can choose.
Clearly, we could reconstruct the whole graph using sufficiently many inputs localized at single nodes, since $Y = A \mathbbm{I}_n$.
However, if we simply aim to identify a cEP with a relatively small number of classes, considerably fewer inputs suffice.

Our idea here is to use input/output behavior as an oracle to simulate the well-known Weisfeiler Lehman algorithm (WL)~\cite{weisfeiler1968reduction}, also known as color refinement.
Starting from an initial coloring $c^{(0)}$ at time $t=0$ (usually the same color for all nodes), the WL algorithm updates the color of each node $v$ iteratively as follows:
$$
c^{(t+1)}(v) = \operatorname{hash}\left(c^{(t)}(v), \{\{ c^{(t)}(x) \mid x \in N(v)\}\}\right)
$$ 
where the doubled brackets denote a ``multi-set'', i.e. a set in which an element can appear more than once. 
Here the hash function is an injective function and ensures that nodes that have (i) the same color in previous iteration $t$ and (ii) the same set of colors in their neighborhood, are assigned the same color in the next iteration $t+1$.
Evidently, every step of the algorithm refines the coloring until at some point the partition of the graph induced by the colors stays the same. 
At this point, all nodes within the same class have the same number of neighbors to each class, i.e., the algorithm found an EP --- the cEP to be precise \cite{paige1987three}.
Using the oracle for $y = Ax$, we present a ``blind'' version of the WL algorithm in~\Cref{alg:BWL}. 
A similar algorithm was proposed by \cite{kersting2014power} for the computation of fractional isomorphisms based on conditional gradients.

\inExtendedVersion{Note that, since we are only interested in nodes with exactly the same colors in their neighborhood, we need only remember if the multi-set is exactly the same or not.  
    A common approach in practice is thus to create an implicit hash function via a dictionary that is indexed by the (sorted) multi-sets. 
Computing the hash of a multi-set then consists of checking whether it is already in the dictionary. If so, one returns the color value stored in the dictionary. If not, the number of entries of the dictionary is stored indexed by the multi-set and returned. This way, distinct multi-sets receive distinct values, but the same multi-set will continue receiving the same value. } 

\begin{algorithm}[t]
  \label{alg:1}
  $O \gets \text{oracle returning } Ax$\;
  $B\gets \mathbbm{1}_{|V|}$; \hfill $\triangleleft$ \textit{start with global partition indicator matrix}\;
  \While{$O(B) \notin \operatorname{span}(B)$} 
  {
	  \For{$i,j \in V$}
	  {
		$B_{i,j} \gets \begin{cases}
			1 & \text{if } O(B)_{i,\_} = O(B)_{j,\_} \wedge B_{i,\_} = B_{j,\_}\\
			0 & \text{else}
		\end{cases}$\
	  }
	remove redundant columns from $B$
  }
  \Return B \hfill $\triangleleft$ \textit{return coarsest equitable partition}
  \caption{BlindWL}
  \label{alg:BWL}
\end{algorithm}

\noindent\textbf{Properties of the blindWL algorithm.}
It is relatively easy to see that~\Cref{alg:BWL} indeed finds the coarsest equitable partition. 
After termination, $AB \in \operatorname{span}(B)$, meaning that there exists some matrix $A^\pi$ with $AB = BA^\pi$.  
Furthermore, consider the cEP represented by $H^*$. 
Throughout the execution, if $B = H^*B^\pi \in \operatorname{span}(H^*)$, then $AB = AH^*B^\pi = H^*A^\pi B^\pi \in \operatorname{span}(H^*)$. 
Since $\mathbbm{1}_{|V|} = H^*\mathbbm{1}_{|H^*|}$, $B$ stays in $\operatorname{span}(H^*)$ and $B$ eventually represents the same EP as $H^*$.

In fact, \Cref{alg:BWL} induces the same partitions as the WL algorithm in each iteration of the while loop: we start with the same color for all nodes, as encoded in the all ones vector.
If the number of neighbors of a certain color $c$ is different for two nodes $u,v$ the WL algorithm will put them in two different classes. 
The corresponding components $(AB)_{u,c} \neq (AB)_{v,c}$ will also be distinct.
Hence, in the next iteration, $B_{u, v} = 0$ meaning $u$ and $v$ are also put into distinct classes by the blind WL algorithm.  

A benefit of the blindWL algorithm is that the intermediate row representations $(AB)_{i,\_}$ of a node $i$ yield an embedding each iteration rather than colors that, if distinct, provide no method of comparison.
For example one can cluster these embeddings to obtain an even coarser partition into nodes that are similar.
This circumvents the sensitivity of the WL algorithm to minor perturbations in the graph: indeed, adding a single edge can disrupt an exact symmetry and yields a much finer cEP. 
While crucial to the original application of the WL (graph isomorphism checking), a more robust approach to assigning the classes is useful for node role extraction and the completely blind problem setting (\cref{sec:experiments}).

\inExtendedVersion{While the proposed algorithm is not as efficient as the actual WL algorithm (which has been thoroughly optimized), it does offer a perspective that can help with detecting other EPs, rather than simply the coarsest EP.
Suppose a graph has multiple distinct EPs, and the cEP is found by the WL algorithm, and the blindWL algorithm starting from the all-ones-vector.  
Note that the blindWL algorithm finds the smallest subspace spanned both by a basis consisting of eigenvectors and a basis consisting of block-standard vectors. 
We know that the all-ones vector must always be part of this subspace, since $H\mathbbm{1}_k = \mathbbm{1}_n$.
When trying to find a finer EP, we need to choose a different set of starting vectors. 
Reasonable candidates can be found by taking an eigenvector that is not a structural eigenvector of the cEP, but that is nonetheless block constant on some nodes, and using the indicator vectors of its blocks.

Assuming that this candidate vector $x$ is in the finer EP, the blindWL algorithm will find the whole EP.}


\section{Scenario II: Truly blind identification of equitable partitions}
\label{sec:main}
We now consider a scenario in which we aim to infer the cEP, but merely observe the outputs of a graph filter excited by a noisy, \emph{random} low-rank excitation over which we have no control.
\begin{equation}
	\label{eq:def_graph_signal}
	y = \alpha f(A)\tilde{H}x + (1-\alpha)z.
\end{equation}
Here $x \sim \mathcal{N}(0,\mathbbm{I}_k)$ and $z \sim \mathcal{N}(0,\mathbbm{I}_n)$ are jointly Gaussian independent random vectors that are each sampled i.i.d from a normal distribution, $\alpha\in[0,1]$ is a parameter that regulates how strongly the structural eigenvectors are excited, and $\tilde{H} = H\operatorname{diag}(1/\sqrt{|C_i|})$. 
For simplicity, we assume that $f(A)$ and $f(A)^2$ have the same cEP indicated by $H$ as $A$. 
Though this seems restrictive, for generic graphs most filters will fulfill this requirement. 
Indeed, the cEP of $A$ is always an EP of $A^k$ as well, though it may not be the coarsest. 

Now observe that the covariance matrix of the above process has the following form:
$$
	\begin{aligned}
		\Sigma &= \mathbb{E}[yy^T] = P \Gamma P^T
		= \alpha^2 f(A) \tilde{H}\tilde{H}^T f(A)^T + (1-\alpha)^2\mathbb{I}_n
	\end{aligned}
$$
Because $f(A)$ has an cEP as indicated by $H$, for any eigenvector $f(A^\pi) v = \lambda v$, associated to a structural eigenvector we have:
\onlyInShortVersion{
$$
	\Sigma Hv = (\alpha^2 \lambda^2 + (1-\alpha)^2) Hv
$$
}\inExtendedVersion{
$$
\begin{aligned}
	\Sigma Hv &= \alpha^2 f(A)\tilde{H}\tilde{H}^Tf(A)Hv + (1-\alpha)^2 \mathbb{I}_n Hv \\
	&= \alpha^2 f(A)\tilde{H}\tilde{H}^T H f(A)^\pi v + (1-\alpha)^2 Hv \\
	&= \alpha^2 f(A) H f(A)^\pi v + (1-\alpha)^2 Hv\\
	&= \alpha^2 H (f(A)^\pi)^2 v + (1-\alpha)^2 Hv\\
	&= (\alpha^2 \lambda^2 + (1-\alpha)^2) Hv,
\end{aligned}
$$
where we have used that $\tilde{H}\tilde{H}^T H = H$.
}Hence, the structural eigenvectors of $\Sigma$ are the same as the structural eigenvectors of $f(A)$, which are scaled-up, block constant eigenvectors of $f(A^\pi)$.
Now, let $f(A) = V\Lambda V^T$ denote the spectral decomposition of the (symmetric) matrix $f(A)$, and denote by $V_\text{EP}$ the subset of structural eigenvectors.

If we consider the cEP $\mathcal{C}^* = \{C_1^*, ..., C_k^*\}$ associated to $f(A)$ and define the $k$-means cost function:
\begin{equation}
	\label{eq:def_F}
	F(\mathcal{C}, V) = \sum_{C \in \mathcal{C}} \sum_{i \in C} \left\lVert V_{i,\_} - \frac{1}{|C|} \sum_{j \in C}V_{j,\_} \right\lVert^2_2
\end{equation}
it is easy to see that $F(\mathcal{C}^*, V_{\text{EP}}) = 0$, as the eigenvectors of $f(A)$ are block-wise constant on the classes of the cEP. 

Hence, if we had access to (a good estimate of) the covariance matrix $\Sigma$, we could simply use $k$-means to find candidates for the cEP, provided we can supply the correct eigenvectors $V$ to the algorithm.
As the above calculations show, the parameter $\alpha$ regulates the scale of the eigenvalues associated to the structural eigenvectors.
For sufficiently large values of $\alpha$ most structural eigenvectors will, in fact, be the dominant eigenvectors.
Assuming that we know the number of classes $k$ of the cEP to be found, we may thus simply pick the top $k$ eigenvectors of $\Sigma$ and optimize the $k$-means objective to obtain the blocks of the EP --- a procedure akin to spectral clustering.
Here we estimate the covariance matrix by the sample covariance based on sampled outputs $y_i$ for $i\in 1,\ldots,s$:
$$
\hat{\Sigma} = \frac{1}{s}\sum_{i=1}^s y_i y_i^T = \hat{P} \hat{\Gamma} \hat{P}^T
$$
For this setup we can show the following result.
\begin{theorem}
	\label{thm:1}
	Let $\hat{\Sigma} - \Sigma = \Delta$ and let $y^1, ..., y^s \in \mathbb{R}^n$ be independent samples from the graph filter as in \cref{eq:def_graph_signal} and let $r= \operatorname{Tr}(\Sigma)/\left\lVert\Sigma\right\lVert_2$.
	Let the following conditions hold:
	\begin{enumerate}
		\item KMeans finds a partition $\hat{\mathcal{C}} = \{\hat{C}_1, ..., \hat{C}_k\}$ that minimizes $F(\hat{\mathcal{C}}, \hat{P}_{\text{EP}})$. 
		\item $\left\lVert y^i\right\lVert_2^2 \leq K \mathbb{E}[\left\lVert y\right\lVert_2^2]$ is bounded almost surely. 
		\item There exists $\delta > 0$ s.t. $\left\lVert \Delta \right\lVert + \delta \leq \gamma_k - \hat{\gamma}_{k+1}$. 
	\end{enumerate}
	Then, for $c > 0$ and with probability at least $1-c$:
		$$
		\sqrt{F(\hat{\mathcal{C}}, V_{\text{EP}})} \leq \frac{\sqrt{8k} \left\lVert \Sigma\right\lVert_2 \Theta \left(\sqrt{\frac{K^2 r \log (n/c)}{s}} + \frac{K^2 r \log (n/c)}{s} \right)}{\delta}
		$$
	for some constant $\Theta$.
\end{theorem}
The proof of the theorem can be found in the full version available \href{https://git.rwth-aachen.de/netsci/blind-extraction-of-equitable-partitions-from-graph-signals}{here}. It is inspired by \cite{wai2018community,wai2019blind} and uses a concentration inequality and the Davis Kahan sin($\theta$) theorem. 
The theorem itself bounds the error of the partition found by the simple spectral clustering algorithm. 
The consistency statement that in the limit $s \rightarrow \infty$ the error vanishes and the extraction of the cEP is exact immediately follows.

A similar statement can be made about an adjusted variant of the blindWL algorithm applied to the estimated covariance matrix $\hat{\Sigma}$.
We can simply replace the exact equality conditions in the computation of the intermediate matrix $B$ in \Cref{alg:BWL} with a clustering algorithm that allows for some variance. 
In the limit $s \rightarrow \infty$, the error of the approximate oracle $\hat{\Sigma}x$ goes toward 0: $\left\lVert\hat{\Sigma} x - \Sigma x\right\lVert = \left\lVert\Delta x\right\lVert \rightarrow 0$. 
Therefore, the adjusted blindWL algorithm also exactly recovers the cEP with infinitely many samples. 
In the next section we explore numerically how the two algorithms compare with finitely many samples.


\section{Experiments}
\label{sec:experiments}

In this section we perform some experiments that support the theoretical findings of this paper. 
Toward this end, we use the setup as described in \cref{sec:main}, \cref{eq:def_graph_signal}. 
While the proposed spectral algorithm of \cref{sec:main} is fit for the task, the blindWL algorithm (\ref{alg:BWL}) must be altered slightly as indicated above. 
In particular, we no longer have control over the inputs, thus the oracle $O$ for $y = Ax$ is replaced with an approximate oracle $\hat{y} = \hat{\Sigma} x$.

Accordingly, we use a (robust) version of the algorithm, in which we replace the exact equality check in line 5 in algorithm \ref{alg:BWL} and instead fit a gaussian mixture on the rows of $O(B)$. 
The adjusted algorithm then uses the indicator vectors of the found clusters as the new intermediate $B$.
We note that in the scenario of \cref{sec:simWL}, both variants of the algorithm yield the same result (under the assumption that the clustering algorithm fits the data optimally).

The graphs used in our synthetic test are sampled from a locally colored configuration model \cite{soderberg2003random}.
As opposed to the original configuration model, in the locally colored configuration model, the edge stubs also specify what color the incident nodes should have.
Specifically, each node in the model has two main sets of parameters: an assignment to a (color) class and a number of colored stubs with which the node is required to link to other classes (which amounts to specifying a partition indicator matrix $H$ and a quotient graph $A^\pi$).
Given that the desired constraints can be met, we obtain a simple graph without self-loops or multi-edges.
Stated differently, using the locally colored configuration model, we can fix the number of colored neighbors for each node and thus guarantee that the sampled graph has an EP $AH = HA^\pi$.
For more details we refer the reader to the code available \href{https://git.rwth-aachen.de/netsci/blind-extraction-of-equitable-partitions-from-graph-signals}{here}. 

In the experiments, graphs with $300$ nodes and an EP with $6$ same sized classes were used. 
In each experiment, we randomly sampled a symmetric matrix $A^\pi$ uniformly from the integers $\{0, ..., 4\}^{6\times 6}$.
Subsequently, we sampled the matrix $A$, generated $s$ outputs $y^i$ (for $i\in 1,\ldots,s$) and evaluated the algorithms.

\begin{figure}[tb]

	\begin{minipage}[b]{1.0\linewidth}
	  \centering
	  \centerline{\includegraphics[width=8.5cm]{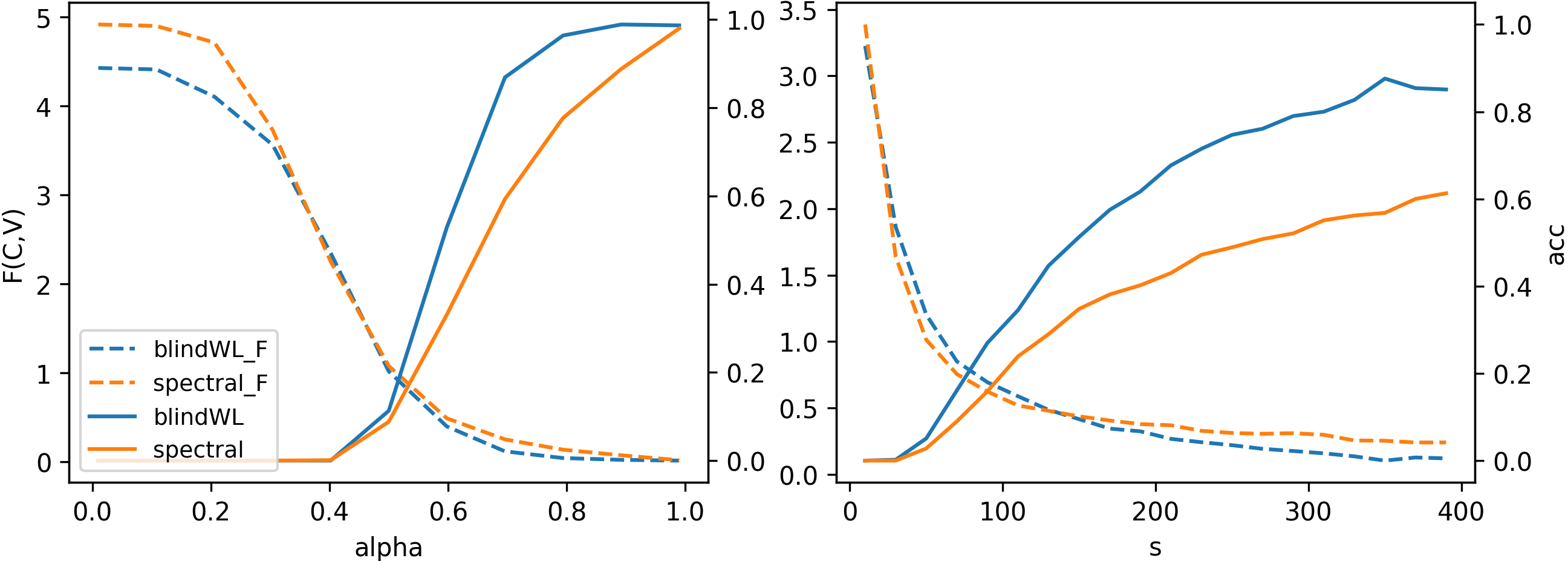}}
	\end{minipage}
	\caption{\textbf{Graph-level and node-level performance of the algorithms.} The figure shows the graph-level accuracy on the right axis (continuous line) and the node-level cost function $F(\hat{\mathcal{C}}, V_{\text{EP}})$ on the left axis (dashed line) of each graph. It also shows the progression of these metrics with fixed sample size $s=300$ and varying noise parameter $\alpha$ (left), and with fixed $\alpha =0.7$ and varying $s$ (right).}
	\label{fig:combined_view}
\end{figure}

We measure the performance of both algorithms using graph-level accuracy, that is, an output partition receives a score of $1$ if it is exactly equivalent to the planted partition; else the score is $0$. 
Note that this is a quite strict measure, as a correct class assignment for all but one node is still counted as a complete failure to recover the EP.
As a second, node-level measure, we use the cost function $F(\hat{\mathcal{C}}, V_{\text{EP}})$ as defined in \cref{eq:def_F}, which can give insight into the quality of wrong partitions.
Both measures are reported as the mean score over $1000$ repeated experiments for each of the parameter configurations shown in \Cref{fig:combined_view}.

In the right plot of~\Cref{fig:combined_view}, a rapid decrease in the cost function and a slightly less steep increase in the accuracy can be seen for increasing sample size, which underlines our theoretical findings in \Cref{thm:1}. 
Though quite close in the node-level measure, the blindWL algorithm already has considerably higher accuracy using only few samples. 

In the left plot of~\Cref{fig:combined_view} and with no noise at all, both algorithms find the correct partitions.
However, the blindWL algorithm is again more robust when increasing the noise. 
The fact that the algorithms do not converge to the same score at $\alpha = 0$ can be explained by the distinct clustering methods. 
While KMeans always finds $6$ clusters the gaussian mixture used in the blindWL can use less than $6$ components in the mixture. 
This should also be kept in mind when comparing the two algorithms, as KMeans requires the number of classes as input, whereas the blindWL algorithm can infer the number of classes from the data.

\section{Conclusion}
\label{sec:conlcusion}
We presented approaches to blindly extracting structural information in the form of an equitable partition of an unobserved graph from only such graph signals.
In \cref{thm:1}, we established a theoretical bound on the error of such an inferred partition and went on to compare the spectral clustering and blindWL approaches experimentally.
An interesting direction for future research may be to exploit this notion of node roles, e.g., the quotient structure may be employed for faster computations of certain graph filters.

\vfill\pagebreak\inExtendedVersion{\vfill\clearpage}

\nocite{boutsidis2015spectral,van1996matrix,vershynin2018high}

\bibliographystyle{IEEEbib}
\bibliography{strings,refs}
\inExtendedVersion{
\vfill\pagebreak
\section{Appendix}
\subsection{Proof of Proposition 1}
\textit{\Cref{prop:1}: If the Perron vector (dominant eigenvector) of a graph is unique, then
it is a structural eigenvector.}
\begin{proof}
	Let $G$ be a graph with adjacency matrix $A$ and its cEP indicated by $H$, that is, $AH = HA^\pi$.
	Consider the dominant eigenpair $(\lambda_1, v_1)$ of $A$. It holds, that:
	$$
		v_1 = \lim_{t \rightarrow \infty} \frac{A^t x}{\lVert A^t x \rVert} 
	$$
	for some $x$, that is not perpendicular to $v_1$. 
	Take $x = \mathbbm{1}_n = H \mathbbm{1}_k$. 
	The all-ones vector $\mathbbm{1}_n$ is not perpendicular to $v_1$, since the dominant eigenvector of a non-negative matrix is non-negative and positive in at least one component.
	Therefore:
	$$
	v_1 = \lim_{t \rightarrow \infty} \frac{A^t \mathbbm{1}_n}{\lVert A^t \mathbbm{1}_n \rVert} =  \lim_{t \rightarrow \infty} \frac{A^t H\mathbbm{1}_k}{\lVert A^t \mathbbm{1}_n \rVert} = \lim_{t \rightarrow \infty} \frac{H(A^\pi)^t \mathbbm{1}_k}{\lVert A^t \mathbbm{1}_n \rVert}
	$$
	Then $v_1 = Hv_1^\pi$ and is thus a structural eigenvector.
	

\end{proof}

\subsection{Proof of Theorem 1}

\textit{
	\Cref{thm:1}. Let $\hat{\Sigma} - \Sigma = \Delta$ and let $y^1, ..., y^s \in \mathbb{R}^n$ be independent samples from the graph filter as in \cref{eq:def_graph_signal} and let $r= \operatorname{Tr}(\Sigma)/\left\lVert\Sigma\right\lVert_2$.
	Let the following conditions hold:
	\begin{enumerate}
		\item KMeans finds a partition $\hat{\mathcal{C}} = \{\hat{C}_1, ..., \hat{C}_k\}$ that minimizes $F(\hat{\mathcal{C}}, \hat{P}_{\text{EP}})$. 
		\item $\left\lVert y^i\right\lVert_2^2 \leq K \mathbb{E}[\left\lVert y\right\lVert_2^2]$ is bounded almost surely. 
		\item There exists $\delta > 0$ s.t. $\left\lVert \Delta \right\lVert + \delta \leq \gamma_k - \hat{\gamma}_{k+1}$. 
	\end{enumerate}
	Then, for $c > 0$ and with probability at least $1-c$:
		$$
		\sqrt{F(\hat{\mathcal{C}}, V_{\text{EP}})} \leq \frac{\sqrt{8k} \left\lVert \Sigma\right\lVert_2 \Theta \left(\sqrt{\frac{K^2 r \log (n/c)}{s}} + \frac{K^2 r \log (n/c)}{s} \right)}{\delta}
		$$
	for some constant $\Theta$.
}

\begin{proof}
	Define the indicator matrices 
	$$\hat{H} = H \operatorname{diag}\left(\left\{\sqrt{|C|} \mid C \in \hat{\mathcal{C}}\right\}\right)$$
	$$ H^* = H \operatorname{diag}\left(\left\{\sqrt{|C|} \mid C \in \mathcal{C}\right\}\right)$$
	It becomes apparent that:
	$$
	\left\lVert\hat{P}_{\text{EP}} - \hat{H}\hat{H}^T\hat{P}_{\text{EP}}\right\lVert_F^2 = F(\hat{\mathcal{C}}, \hat{P}_{\text{EP}})
	$$
	and by condition 1, $\hat{H}$ minimizes the expression. 
	Similarly, 
	$$
	\left\lVert V_{\text{EP}} - H^*(H^*)^T V_{\text{EP}}\right\lVert_F^2 = F(\mathcal{C}^*, V_{\text{EP}}) = 0
	$$
	In the following the aim will be to bound 
	$$
	\sqrt{F(\hat{\mathcal{C}}, V_{\text{EP}})} = \left\lVert V_{\text{EP}} - \hat{H}\hat{H}^T V_{\text{EP}}\right\lVert_F
	$$
	As a shorthand, define the error matrix $E = V_{\text{EP}}V_{\text{EP}}^T - \hat{P}_{\text{EP}}\hat{P}_{\text{EP}}^T$. The following inequalities hold:
	$$
	\begin{aligned}
	&\left\|V_{\text{EP}}-\hat{H} \hat{H}^{\top} V_{\text{EP}}\right\|_{\mathrm{F}}\\
	&=\left\|\left(I-\hat{H} \hat{H}^{\top}\right) V_{\text{EP}} \right\|_{\mathrm{F}}\\
	&=\left\|\left(I-\hat{H} \hat{H}^{\top}\right) V_{\text{EP}} V_{\text{EP}}^{\top}\right\|_{\mathrm{F}}\\
	&=\left\|\left(I-\hat{H} \hat{H}^{\top}\right)\left(\hat{P}_{\text{EP}} \hat{P}_{\text{EP}}^{\top}+E\right)\right\|_{\mathrm{F}} \\
	&\leq\left\|\left(I-\hat{H} \hat{H}^{\top}\right) \hat{P}_{\text{EP}} \hat{P}_{\text{EP}}^{\top}\right\|_{\mathrm{F}}+\left\|\left(I-\hat{H} \hat{H}^{\top}\right) E\right\|_{\mathrm{F}} \\
	&\leq\left\|\left(I-\hat{H} \hat{H}^{\top}\right) \hat{P}_{\text{EP}} \hat{P}_{\text{EP}}^{\top}\right\|_{\mathrm{F}}+\|E\|_{\mathrm{F}}
	\end{aligned}
	$$
	where the second equality stems from the fact that $V$ is unitary and the frobenius norm is invariant under unitary operations. 
	In the same vein, the last equality is a consequence of $(I - \hat{H}\hat{H}^T)$ being a projection matrix and therefore also having no influence on the norm.
	Now, since $\hat{\mathcal{C}}$ minimizes $F(C, \hat{P}_{\text{EP}})$ going to the ground-truth partition will actually increase the cost:
	$$
	\begin{aligned}
	&\left\|\left(I-\hat{H} \hat{H}^{\top}\right) \hat{P}_{\text{EP}} \hat{P}_{\text{EP}}^{\top}\right\|_{\mathrm{F}}+\|E\|_{\mathrm{F}} \\
	&\leq\left\|\left(I-H^{*}\left(H^{*}\right)^{\top}\right) \hat{P}_{\text{EP}} \hat{P}_{\text{EP}}^{\top}\right\|_{\mathrm{F}}+\|E\|_{\mathrm{F}} \\
	&=\left\|\left(I-H^{*}\left(H^{*}\right)^{\top}\right)\left(V_{\text{EP}} V_{\text{EP}}^{\top}-E\right)\right\|_{\mathrm{F}}+\|E\|_{\mathrm{F}} \\
	&\leq\left\|\left(I-H^{*}\left(H^{*}\right)^{\top}\right) V_{\text{EP}} V_{\text{EP}}^{\top}\right\|_{\mathrm{F}}+2\|E\|_{\mathrm{F}} \\
	&=\sqrt{F(\mathcal{C}^*, V_{\text{EP}})}+2\|E\|_{\mathrm{F}}\\
	&= 2\|E\|_{\mathrm{F}}
	\end{aligned}
	$$
	To now bound the error $\left\lVert E\right\lVert_F$ we change the norms by the following lemma:
	\begin{lemma}[\cite{boutsidis2015spectral}, Lemma 7]
		For any $A, B \in \mathbb{R}^{m \times n}$ with $m \geq n$ and $A^TA = B^TB = I_n$, it holds that:
		$$
		\left\lVert AA^T - BB^T\right\lVert_F^2 \leq 2 n \left\lVert AA^T - BB^T\right\lVert_2^2
		$$
	\end{lemma}
	\noindent It directly follows that:
	$$
	\sqrt{F\left(\hat{\mathcal{C}}, V\right)} \leq 2\sqrt{2k} \left\lVert E\right\lVert_2
	$$
	We now need to bound the error term $E = V_{\text{EP}}V_{\text{EP}}^T - \hat{P}_{\text{EP}}\hat{P}_{\text{EP}}^T$. 
	Since the first term consists of the eigenvectors of the cEP of $A$, which are by assumption the same as the eigenvectors of the cEP of $A^2$, which are, in turn, the same as the eigenvectors of $\Sigma$, we have $E = P_{\text{EP}}P_{\text{EP}}^T - \hat{P}_{\text{EP}}\hat{P}_{\text{EP}}^T$.
	As a direct consequence of \cite[Theorem 2.6.1]{van1996matrix} this, in turn, can be rewritten as: 
	$$
	\left\lVert E\right\lVert_2 = \left\lVert P_{\text{EP}}P_{\text{EP}}^T - \hat{P}_{\text{EP}}\hat{P}_{\text{EP}}^T\right\lVert_2 = \left\lVert\hat{P}_{\overline{\text{EP}}}^TP_{\text{EP}}\right\lVert_2
	$$
    We continue the proof by applying a variant of the Davis Kahan sin $\theta$ theorem to show that:
	$$
	\left\lVert \hat{P}_{\overline{\text{EP}}}^TP_{\text{EP}}\right\lVert_2 \leq  \frac{\left\lVert \hat{P}_{\overline{\text{EP}}}^T\Delta P_{\text{EP}}\right\lVert_2}{\delta} \leq \frac{\left\lVert\Delta\right\lVert_2}{\delta}
	$$
	The second inequality is obvious as the matrices except $\Delta$ are orthogonal and therefore do not increase the largest eigenvalue. 
	For the first inequality, consider the decomposition of the covariance matrix:
	$$
	\Sigma = P \Gamma P^T = P_{\text{EP}} \Gamma_{\text{EP}} P_{\text{EP}}^T + P_{\overline{\text{EP}}} \Gamma_{\overline{\text{EP}}} P_{\overline{\text{EP}}}^T 
	$$ 
	The estimated covariance matrix is also decomposed in the same way. Now, it holds that:
	$$
	 \begin{aligned}
		\Delta P_{\text{EP}} &= \left(\Sigma + \Delta\right) P_{\text{EP}} - \Sigma P_{\text{EP}}\\
		 &= \hat{\Sigma} P_{\text{EP}} - \left(P_{\text{EP}} \Gamma_{\text{EP}} P_{\text{EP}}^T + P_{\overline{\text{EP}}} \Gamma_{\overline{\text{EP}}} P_{\overline{\text{EP}}}^T\right) P_{\text{EP}}\\
		 &= \hat{\Sigma} P_{\text{EP}} - P_{\text{EP}} \Gamma_{\text{EP}}
	 \end{aligned}
	$$
	Furthermore:
	$$
	\begin{aligned}
	\hat{P}_{\overline{\text{EP}}}^T \Delta P_{\text{EP}} &= \hat{P}_{\overline{\text{EP}}}^T \hat{\Sigma} P_{\text{EP}} - \hat{P}_{\overline{\text{EP}}}^T P_{\text{EP}} \Gamma_{\text{EP}}\\
	&= \hat{\Gamma}_{\overline{\text{EP}}} \hat{P}_{\overline{\text{EP}}}^T P_{\text{EP}} - \hat{P}_{\overline{\text{EP}}}^T P_{\text{EP}} \Gamma_{\text{EP}}\\
	\end{aligned}
	$$
	We will now use a trick to center the eigenvalues of $\Gamma$ around 0:
	Let $c = (\gamma_1 + \gamma_{k})/2$ and $d = (\gamma_1 - \gamma_k)/2$. By the triangle inequality:
	$$
	\begin{aligned}
		\left\lVert P_{\overline{\text{EP}}}^T \Delta P_{\text{EP}}\right\lVert_2 &= \left\lVert\hat{\Gamma}_{\overline{\text{EP}}} \hat{P}_{\overline{\text{EP}}}^T P_{\text{EP}} - \hat{P}_{\overline{\text{EP}}}^T P_{\text{EP}} \Gamma_{\text{EP}}\right\lVert_2\\
		&= \Big\lVert\left(\hat{\Gamma}_{\overline{\text{EP}}} - c\mathbb{I}_n\right) \hat{P}_{\overline{\text{EP}}}^T P_{\text{EP}} + c\mathbb{I} \hat{P}_{\overline{\text{EP}}}^T P_{\text{EP}} \\
		& \quad - \hat{P}_{\overline{\text{EP}}}^T P_{\text{EP}} \Big(\Gamma_{\text{EP}} - c\mathbb{I}_n\Big) - \hat{P}_{\overline{\text{EP}}}^T P_{\text{EP}} c\mathbb{I}_n \Big\lVert_2\\
		&= \Big\lVert\left(\hat{\Gamma}_{\overline{\text{EP}}} - c\mathbb{I}_n\right) \hat{P}_{\overline{\text{EP}}}^T P_{\text{EP}}\\
		&\quad\quad - \hat{P}_{\overline{\text{EP}}}^T P_{\text{EP}} \Big(\Gamma_{\text{EP}} - c\mathbb{I}_n\Big)\Big\lVert_2\\
		&\geq \left\lVert\left(\hat{\Gamma}_{\overline{\text{EP}}} - c\mathbb{I}_n\right) \hat{P}_{\overline{\text{EP}}}^T P_{\text{EP}}\right\lVert_2\\
		&\quad\quad- \left\lVert\hat{P}_{\overline{\text{EP}}}^T P_{\text{EP}} \Big(\Gamma_{\text{EP}} - c\mathbb{I}_n\Big)\right\lVert_2\\
	\end{aligned}
	$$
	As $\Gamma$ is a diagonal matrix and $\gamma_i - c \in [-d,d]$ for $i \in \{1, ..., k\}$, it holds that $\left\lVert\Gamma_{\text{EP}} - c \mathbb{I}_n\right\lVert_2 \leq d$. 
	Applying this yields:
	$$
	\left\lVert\hat{P}_{\overline{\text{EP}}}^T P_{\text{EP}} \Big(\Gamma_{\text{EP}} - c\mathbb{I}_n\Big)\right\lVert_2 \leq d \left\lVert\hat{P}_{\overline{\text{EP}}}^T P_{\text{EP}}\right\lVert_2 
	$$
	Using the same logic, and by condition 3, $\hat{\gamma}_i - c \in [-\infty, -(d+\delta)]$ for $i \in \{k+1, ..., n\}$. Thus, $\left\lVert\hat{\Gamma}_{\overline{\text{EP}}} - c \mathbb{I}_n\right\lVert_2 \geq (d+\delta)$.
	Applying this to the norms, it holds that:
	$$
	\left\lVert \Big(\hat{\Gamma}_{\overline{\text{EP}}} - c\mathbb{I}_n\Big) \hat{P}_{\overline{\text{EP}}}^T P_{\text{EP}}\right\lVert_2 \geq (d+\delta) \left\lVert\hat{P}_{\overline{\text{EP}}}^T P_{\text{EP}}\right\lVert_2
	$$
	Subtracting the right-hand side of the equations from one another, we end up with the inequality:
	$$
	\begin{aligned}
	\left\lVert P_{\overline{\text{EP}}}^T \Delta P_{\text{EP}}\right\rVert_2 &\geq (d+\delta) \left\lVert\hat{P}_{\overline{\text{EP}}}^T P_{\text{EP}}\right\rVert_2 - d \left\lVert\hat{P}_{\overline{\text{EP}}}^T P_{\text{EP}}\right\rVert_2 \\
	&\geq \delta \left\lVert\hat{P}_{\overline{\text{EP}}}^T P_{\text{EP}}\right\rVert_2
	\end{aligned}	
	$$
	It has now been proven, that
	$$
	\sqrt{F(\hat{\mathcal{C}}, V)} \leq 2 \sqrt{2k} \frac{\left\lVert\Delta\right\lVert_2}{\delta}
	$$
	Applying the following theorem then directly yields the theorem stated above.
	\begin{theorem}[\cite{vershynin2018high} Theorem 5.6.1, 5.6.4]
		Let $y^1, ..., y^s$ be independent samples of the graph filter as in \cref{eq:def_graph_signal}, let $\left\lVert y^l\right\lVert_2 \leq K \mathbb{E}[\left\lVert y\right\lVert_2^2]$ and let the effective rank of the covariance be $r = \operatorname{Tr}(\Sigma)/\left\lVert\Sigma\right\lVert_2$, then for every $c>0$ it holds, that:
		$$
			\left\lVert\Delta\right\rVert_2 \leq \left\lVert\Sigma\right\rVert_2 C \left(\sqrt{\frac{K^2 r \log (n/c)}{s}} + \frac{K^2 r \log (n/c)}{s} \right)
		$$
		with probability at least $1-c$.
	\end{theorem}
\end{proof}
}

\end{document}